\documentclass{sn-jnl}
\usepackage[numbers,sort&compress]{natbib}

\usepackage{graphicx}%
\usepackage{multirow}%
\usepackage{amsmath,amssymb,amsfonts}%
\usepackage{amsthm}%
\usepackage{mathrsfs}%
\usepackage[title]{appendix}%
\usepackage{xcolor}%
\usepackage{textcomp}%
\usepackage{manyfoot}%
\usepackage{booktabs}%
\usepackage{siunitx}
\usepackage[noend]{algpseudocode}%
\usepackage{listings}%
\usepackage{placeins}
\usepackage[linesnumbered,noend]{algorithm2e}

\SetCommentSty{mycommfont}
\RestyleAlgo{ruled}
\SetKwComment{Comment}{/* }{ */}
\usepackage{hyperref}
\usepackage[nameinlink]{cleveref} %

\usepackage{lmodern} %

\usepackage{tabularx}
\usepackage{orcidlink}

\newcommand{\T}{\mathsf{\scriptscriptstyle T}}
\newcommand*{\defeq}{\mathrel{\vcenter{\baselineskip0.5ex \lineskiplimit0pt
                     \hbox{\footnotesize.}\hbox{\footnotesize.}}}%
                     =}

\usepackage{bbm}
\newcommand{\field}[1]{\mathbbm{#1}}
\newcommand{\reals}{\field{R}}
\newcommand{\integers}{\field{Z}}

\newcommand{\nonnegreals}{\reals_{\scriptscriptstyle \ge 0}}
\newcommand{\nonnegints}{\integers_{\scriptscriptstyle \ge 0}}
\newcommand{\nonnegintegers}{\nonnegints}

\newcommand{\R}{\reals}

\newcommand{\Z}{\integers}

\newcommand{\figurewidth}{\textwidth}
\newcommand{\instance}[1]{\texttt{#1}}
\newcommand{\latin}[1]{\textit{#1}}

\begin{document}

\title[Machine Learning for Optimization-Based Separation of Mixed-Integer Rounding Cuts]{Machine Learning for Optimization-Based Separation of Mixed-Integer Rounding Cuts}

\author*[1]{\fnm{Oscar} \sur{Guaje}\,\orcidlink{0000-0002-4110-6958}}\email{o.guaje@mail.utoronto.ca}

\author[1]{\fnm{Arnaud} \sur{Deza}\,\orcidlink{0009-0008-2743-9298}}\email{arnaud.deza@mail.utoronto.ca}

\author*[2]{\fnm{Aleksandr M.} \sur{Kazachkov}\,\orcidlink{0000-0002-4949-9565}}\email{akazachkov@ufl.edu}

\author[1]{\fnm{Elias B.} \sur{Khalil}\,\orcidlink{0000-0001-5844-9642}}\email{khalil@mie.utoronto.ca}

\affil[1]{\orgdiv{Department of Mechanical \& Industrial Engineering}, \orgname{University of Toronto}, \orgaddress{\city{Toronto}, \state{ON}, \country{Canada}}}

\affil[2]{\orgdiv{Department of Industrial \& Systems Engineering}, \orgname{University of Florida}, \orgaddress{\city{Gainesville}, \state{FL}, \country{USA}}}

\abstract{%
Mixed-integer rounding (MIR) cutting planes (cuts) are effective at improving the strength of a linear relaxation for mixed-integer linear programming (MIP) problems. The cuts in this family are derived by aggregating constraints then rounding coefficients, but finding the strongest MIR cuts requires optimizing a costly MIP for the aggregation step, so in practice, heuristic strategies for separating fractional points are employed. We propose to improve MIR cut generation in the context of a common scenario in applications, where constraints remain fixed but costs are varied. We present a hybrid cut generation framework in which we train a machine learning (ML) model to classify which constraints are involved in useful MIR cuts based on fractional points from relaxations of the problem. At test time, the predictions of the ML model create a reduced MIP-based generator of MIR cuts. In our experiments, we create an instance family from each of three benchmark MIP instances by performing a careful and costly perturbation of objective coefficients to build a dataset of 1,000 fractional points to be separated over the same constraint set. The results indicate that the reduced separator better strengthens the bound in each round of cut generation, particularly for instances in which the full separator failed to find strong cuts.
}


\maketitle

\section{Introduction}

Cutting planes, or cuts, are used by mixed-integer linear programming (MIP) solvers to strengthen formulations and improve running time. Some cut families contain a large number of possible inequalities, leading to heuristic separation strategies and prompting research into optimization-based generation to find the best cuts~\cite{cgmip,procgmip,Balas2008,MIRseparator}. In this work, we investigate a learning-based approach to reduce the computational cost of MIP subproblems for generation of \emph{mixed-integer rounding} (MIR) inequalities~\cite{wolsey1999integer,MIR}.

An MIR cut is derived via a weighted aggregation of a subset of constraints, followed by appropriate rounding of the resulting base inequality's coefficients.
Selecting which constraints to aggregate is, in practice, based on computationally-inexpensive heuristics~\cite{Tobias,Lodi2013},
which may miss strong MIR cuts obtainable by optimization.
However, applications often involve repeatedly solving related optimization problems, which provides an opportunity for inferring helpful structural properties.
Specifically, we consider a practitioner that solves a set of MIP instances with the same constraints but varying cost vectors.

We ask the following question: assuming access to samples of a distribution of similar MIP instances, can one \textit{learn} useful constraint selection models from an expensive optimization-based generator for MIR cuts? 
Our work tests two hypotheses:
    (i) information about effective cutting planes for a distribution of MIP instances can inform the separation of cutting planes for similar but previously unseen MIP instances; and 
    (ii) an optimization-based generator can be accelerated by carefully fixing to zero some of its decision variables, i.e., eliminating some constraints from consideration, at little sacrifice to the dual bound improvement obtained by the original separator.
Our approach can be seen as a hybrid of machine learning (ML) and MIP: an ML model accelerates a MIP-based generator.

We frame the ML problem as supervised binary classification of constraints. The MIP-based separator of~\citet*{MIRseparator} is executed for a number of rounds on each training instance to produce a portfolio of cuts, tracking which constraints are aggregated to produce each cut. 
Binary ``labels'' are assigned to each constraint depending on whether it was selected to generate a good MIR cut or not. Given a MIP instance, a fractional vertex to be separated, and a particular constraint, we devise and compute a set of 54 features that contextualize the relationship between that constraint and the vertex, as well as that constraint and other constraints of the problem. The supervised learning problem is then one of finding an accurate mapping from features to labels. Given an unseen test instance in which a different fractional point is to be cut, the trained ML classifier selects a subset of promising constraints, after which MIP-based MIR separation is restricted only to those constraints.

We test our framework on datasets derived from three benchmark instances from the \emph{2017 Mixed Integer Programming Library} (MIPLIB 2017)~\cite{miplib} that have nontrivial \emph{integrality gaps},
i.e., the objective value of an integer-optimal solution is significantly different than the optimal value after relaxing integrality constraints.
We create a large dataset from each base MIPLIB 2017 instance by perturbing the objective function coefficients to find {1,000} different initial fractional solutions that need to be cut off.
Our results indicate that learning a useful classifier is indeed possible and can result in favorable performance compared to running the separator with all constraints.

\section{Related Work}
In recent years, several papers have explored ML techniques to enhance optimization solvers. In the context of general-purpose MIP solvers, ML has been used to design branching strategies \cite{Elias} and to decide when and which heuristics to run \cite{EliasHeuristics}. As for cutting planes, most work has been on selecting which cuts to add from a pool of cuts \cite{HUANG2022,columbia,turner2022}, while a more recent paper has tackled the problem of deciding whether to generate cuts only at the root node of the branch-and-bound tree~\cite{berthold2022}. 

One stream of literature has studied Chv{\'{a}}tal-Gomory (CG) cuts, which also involve aggregation of constraints followed by a different rounding step~\cite{CHVATAL}.
\citet{NEURIPS2021_210b7ec7} study the sample complexity of learning CG constraint aggregation coefficients. 
They establish lower bounds on the number of samples (instances) needed to accurately estimate expected search tree size resulting from adding CG cuts with specific coefficients.
\citet{deza2024} present a framework to predict a set of useful CG coefficients for aggregation using graph neural networks. \citet{becu2024} use the aggregation weights of known instance variations to Gomory mixed integer cuts for other variations of the same instance.
\citet{dragotto2023} generate more general split cuts with the help of ML, identifying a split disjunction by a neural network.

A more comprehensive survey of ML for cuts literature was recently completed by \citet{deza2023machine}. We note that ML-based approaches have also been used to exploit problem structures and accelerate the solving process for specific classes of problems (see, e.g., \citet{Larsen,unitcomm}).

\section{Methodology}
\label{sec:methodology}

In this section, we describe how to optimize over MIR cuts, add in a step to produce a pool of cuts, introduce our data generation procedure, set up our classification model, and explain how this is then given to a reduced MIR generator.

We consider a MIP problem over integer variables $x \in \nonnegintegers^n$ and continuous variables $v \in \nonnegreals^p$ associated to rational objective coefficients $f$ and $g$, subject to $m$ constraints with rational data $A$, $C$, $b$:

\begin{equation}
    \label{eq:mip}
    \begin{aligned}
    \min_{x,v} 
            \quad &f^\T x + g^\T  v
            \\
            &Ax + Cv = b,
            \\
            &x,v \ge 0,
            \\
            &x\in \Z^n,
            \\
            &v \in \R^p
            .
    \end{aligned}
\end{equation}

Let $(x^\star, v^\star)$ be a fractional feasible solution to the linear relaxation of \eqref{eq:mip}, obtained by relaxing the integrality restrictions.
Our target is to find cuts that are violated by $(x^{\star}, v^{\star})$ but valid for \eqref{eq:mip}, i.e., satisfied by all integer-feasible points.

\subsection{Generation of MIR Cuts}

We first provide an overview of MIR inequalities, following the presentation by \citet*{MIRseparator}, based on \cite{NemWol90_recursive-procedure} and \cite{Wolsey98_integer-programming}.
We then review the optimization-based approach by \citet*{MIRseparator} to identifying strong cuts in the family.

\subsubsection{Relaxed MIR Inequalities}

Let $\lambda \in \reals^m$ be a set of multipliers that we treat as a row vector to aggregate the constraints of the MIP \eqref{eq:mip} into a single constraint
\begin{equation}
    \lambda A x + \lambda C v = \lambda b.
    \label{eq:agg-constraint}
\end{equation}

Suppose $\hat{\alpha} \in [0,1]^n$ and $\bar{\alpha} \in \Z^n$ satisfy
    $\hat{\alpha} + \bar{\alpha} \ge \lambda A$,
where $\hat{\alpha}$ represents the fractional component.
Further, let $c^{+} \in \nonnegreals^{p}$ with $c^{+} \ge \max \{ 0, \lambda C \}$.
Finally, let $\hat{\beta} \in [0,1]$ and $\bar{\beta} \in \Z$ such that
    $\hat{\beta} + \bar{\beta} \le \lambda b$,
in which $\hat{\beta}$ is again the fractional component.
Then, for any $(x,v)$ feasible to the linear relaxation of \eqref{eq:mip}, it holds that
\begin{equation*}
        (\hat{\alpha} + \bar{\alpha}) x
        + c^{+} v
    \ge
        \lambda A x + \lambda C v
    =
        \lambda b
    \ge
        \hat{\beta} + \bar{\beta}.
\end{equation*}
We can then deduce the validity of the following ``relaxed'' MIR inequality~\cite[Eq.~(11)]{MIRseparator}:
    \begin{equation*}
        \hat{\alpha} x +
        \hat{\beta} \Bar{\alpha} x + 
        c^{+}v
        \ge \hat{\beta}\left(\Bar{\beta} + 1 \right).
    \end{equation*}

\subsubsection{Separating from the Family of MIR Inequalities}
MIR cuts are effective in practice at strengthening a MIP's linear relaxation~\cite{MIR,Tobias}.
However, as there are infinitely many aggregation vectors $\lambda$, each of which yields a valid cut, finding the ``best'' cut to add is a nontrivial task.
One approach to computing a ``good'' MIR cut is to formulate a MIP to find an aggregation vector that maximizes some measure of violation of a cut by a fractional solution. Given a MIP instance and a fractional solution $(x^\star,v^\star)$, \citet*[Eqs.~(12)--(18)]{MIRseparator} propose the following MIP, referred to as Appx-MIR-Sep, to (approximately) search for the MIR cut most violated by $(x^\star, v^\star)$:
\begin{subequations}
	\label{eq:mir_sep}
	\begin{align}
		\max_{\lambda,\hat{\alpha},\bar{\alpha},c^{+},\hat{\beta},\bar{\beta},\{\Delta_k,\pi_k\}_{k \in K},\Delta}
        \quad 
            &\sum_{k\in K} \varepsilon_k \Delta_k  - \left(c^{+}v^{\star} + \hat{\alpha} x^{\star}\right) \label{eq:sep_obj}                         \\
		&\hat{\alpha} + \Bar{\alpha} \geq \lambda A                                                                                      \label{eq:sep_alpha}                       \\
		&c^{+}    \geq \lambda C                                                                                     \label{eq:sep_cpos}                        \\
	    &\hat{\beta} + \Bar{\beta}   \leq \lambda b                                                                                     \label{eq:sep_beta}                        \\
		&\hat{\beta}                 \geq \sum_{k\in K} \varepsilon_k \pi_k                                                             \label{eq:sep_lin_beta}                    \\
		&\Delta                      = \left( \Bar{\beta} + 1 \right) - \Bar{\alpha} x^{\star}                                          \label{eq:sep_def_delta}                   \\
		&\Delta_k                    \leq \Delta                                                                                       &&\forall \ k \in K \label{eq:sep_lin_delta} \\
		&\Delta_k                    \leq \pi_k                                                                                         &&\forall \ k \in K \label{eq:sep_lin_pi}    \\
		&c^{+} \ge 0 \label{eq:sep_dom_c} \\
        &\hat{\alpha}, \hat{\beta} \in [0,1]
        \label{eq:sep_dom_hatalpha_hatbeta}                   \\
		&\pi                         \in \{ 0,1 \}^{|K|}                                                                                \label{eq:sep_dom_pi}                      \\
		&\Bar{\alpha},\Bar{\beta}                 \in \mathbb{Z}.                                                                                    \label{eq:sep_int_alpha_beta}
	\end{align}
\end{subequations}

Since computing the violation of an MIR cut would yield a nonlinear model, Appx-MIR-Sep approximates the violation $\varepsilon$ of a cut by assuming it is representable over a set $\mathcal{E}=\lbrace \varepsilon_k=2^{-k}: k\in K \rbrace$. Any feasible solution to Appx-MIR-Sep can be used to compute an MIR cut to add to the linear relaxation of the problem. The reader is referred to \cite{MIRseparator} for a more complete account of this model. 

What is of interest to us here is that MIR cuts empirically close large gaps and can be separated via a MIP. 
As cut separation is most useful at the root node of a branch-and-bound tree, the \textit{de facto} procedure for MIP solving, it is rather impractical to solve another complex separation MIP to facilitate solving the original MIP instance. However, since the number of variables in Appx-MIR-Sep depends on the number of constraints in the original MIP, a straightforward way to reduce Appx-MIR-Sep's computational cost is to reduce the  number of constraints to aggregate (i.e., fix some entries of the vector $\lambda$ to zero \latin{a priori}). 
This is what we will attempt to do by training an ML classifier that identifies ``useful'' versus ``unimportant'' constraints.

\subsection{Populating a Pool of Cuts} \label{sec:pool}

Given a fractional point, Appx-MIR-Sep seeks the most violated MIR cut. However, violation is one of many cut quality metrics~\cite{Wess}.
An arguably more relevant metric, for which violation serves as a tractable but imperfect substitute, is \emph{percent (integrality) gap closed}, defined as the percent of the ``gap'' between the optimal objective value of the MIP ($z_{I}$) and of its LP relaxation ($z_{LP}$) that is closed by adding cuts (leading to objective value $z$):
    $\text{GapClosed} \defeq 100 \cdot {(z - z_{LP})} / {(z_I - z_{LP})} $.

To increase the chance of obtaining MIR cuts with large gap closed via \mbox{Appx-MIR-Sep}, we use the fact that any feasible solution to Appx-MIR-Sep gives a valid cut, and populate a ``pool'' of cuts using the off-the-shelf capability of modern solvers to collect \textit{multiple} feasible solutions to a MIP. These alternative (and potentially suboptimal for \mbox{Appx-MIR-Sep}) cuts may result in a larger gap closed than the most violated cut. Similarly to the approach by \citet*{procgmip} and \citet*{MIRseparator}, we populate the cut pool with all the incumbent solutions that the solver finds while solving Appx-MIR-Sep within a given time limit, add all the cuts in the pool to the linear relaxation, get a new fractional solution to separate, and repeat this procedure iteratively until Appx-MIR-Sep is unable to find a separating cut. \Cref{alg:full} describes the complete cutting loop; \cref{step:ml:get-features,step:ml:get-constraints} are skipped when the classification model is not used, resulting in the benchmark \textbf{``full'' separator} that operates on all constraints.

\begin{algorithm}[hbt!] 
\caption{Cutting loop (optionally with ML)}\label{alg:full}
\SetKwInOut{Input}{Input}
\Input{
    MIP instance $P$ in the form~\eqref{eq:mip};
    ML model $\zeta$ (optional)
}
Initialize: {$LP \gets \text{LP relaxation of $P$}$}; {$(x^\star,v^\star) \gets \text{Solve}(LP)$}

\While{$x^\star \not\in \mathbb{Z}^n $}{
  $\texttt{features}_\mathcal{C} \gets \texttt{get\_features}(P, x^\star, v^\star)$ \label{step:ml:get-features}
  
  $\Lambda \gets \zeta(\texttt{features}_\mathcal{C})$ \Comment*[r]{$\Lambda$ are constraints $\zeta$ predicts to be useful} \label{step:ml:get-constraints}
  
  $\mathcal{C} \gets \text{Solve Appx-MIR-Sep}(x^\star, v^\star, \Lambda)$
  \Comment*[r]{$\mathcal{C}$ is a set of cuts}

  \eIf{$\mathcal{C}=\emptyset $}
    {\textbf{break}}
    {
        {$LP \gets LP \cup \mathcal{C}$  
        \Comment*[r]{Add cuts to LP}
    }
    {$(\hat{x}, \hat{v}) \gets \text{Solve}(LP)$}
 
    {\eIf{($\hat{x} , \hat{v}) = (x^\star, v^\star)$}
    {\textbf{break}}
    {$(x^\star, v^\star) \gets (\hat{x}, \hat{v})$}}
  }
}
\end{algorithm}

\subsection{Learning for MIR Cuts}

\subsubsection{Data Generation}

To conduct supervised learning, one must assume access to a sample over a distribution of similar MIP instances. While these may be available in a variety of application domains or through synthetic generators for specific families of combinatorial problems, we have opted to generate perturbations of instances from the MIPLIB 2017 library of benchmark MIP instances~\cite{miplib}. 
We will now detail that generation process.

Consider a ``base'' MIP instance $P$ in the form~\eqref{eq:mip}.
If the base MIP instance contains inequality constraints, we explicitly add continuous slack variables, and we add variable upper bounds as rows of the constraint matrix.
We generate a new instance 
    $\tilde{P} \defeq \min \lbrace  \tilde{f}^\T x +  \tilde{g}^\T  v:  {A}x + {C}v =  {b}, x\in \mathbb{Z}^n, v\in\mathbb{R}^p, x,v\geq0 \rbrace$
by replacing the objective coefficients with random vectors $\tilde{f}$ and $\tilde{g}$ drawn from a truncated normal distribution,
then using rejection sampling to ensure that $\tilde{P}$ has a different optimal solution to its linear relaxation compared to $P$.
Specifically, every positive component of $(f,g)$ is replaced by $\max\{0,u^{+}\}$, where $u^{+}$ is normally distributed with mean and variance computed based on all positive components of $(f,g)$,
while every negative component of $(f,g)$ is replaced by $\min\{0,u^{-}\}$ where $u^{-}$ is normally distributed with mean and variance computed based on the negative components of $(f,g)$.

We denote by $\mathcal{P}$ the \textit{instance family} generated from the base instance $P$, where each member of the family is identified by a different objective vector.
For any pair of instances from the same family, Appx-MIR-Sep only differs in the objective function~\eqref{eq:sep_obj} and in constraint~\eqref{eq:sep_def_delta}. 
Following \Cref{alg:full}, for every instance $P\in \mathcal{P}$, we find a solution $(x^\star, v^\star)^0$ to its linear relaxation. 
If this solution is not integer-feasible, we solve \mbox{Appx-MIR-Sep} to obtain a set of MIR cuts that we add to the linear relaxation, for which we then compute a new solution $(x^\star, v^\star)^1$. 
We repeat these \emph{rounds of cuts} until we find a feasible solution to the original MIP instance, or Appx-MIR-Sep cannot find a cut that separates the last fractional solution. We generate only \emph{rank-1} cuts, i.e., previously-computed cuts are not used in the derivation of new cuts, so the only difference between Appx-MIR-Sep in different rounds is the objective function~\eqref{eq:sep_obj} and constraint~\eqref{eq:sep_def_delta}. Finally, after each run of Appx-MIR-Sep, we record the multiplier vector $\lambda$ for each cut in the solution pool.

\subsubsection{Classification Models}

To test both our hypotheses, we train an ML model that uses information about a MIP instance and a fractional solution and predicts a subset of constraints that will be useful in generating MIR cuts. The prediction is made for each constraint of the instance independently.

Each observation in our dataset corresponds to a constraint of an instance $\tilde{P}$ and a round of separation.
For example, for an instance family $\mathcal{P}$ with $\lvert \mathcal{P} \rvert =10$, constructed from a base instance $P$ with $m=5$ constraints, for which the cutting loop did 2 rounds, the dataset would have $10 \times 5 \times 2=100$ observations. We construct a set of features based on previous work on ML for MIP \cite{Elias}, and on measures that are traditionally used to score cuts \cite{Wess}, which we describe in~\Cref{tab:features}. These features describe the instance (e.g., through statistics of the cost coefficients), the constraint (e.g., the constraint's constant side $b_i$, $1 \le i \le m$), and the relationship between the fractional point of interest and the constraint (e.g., the value of the corresponding slack variable, the distance between the constraint's hyperplane and the point). This results in a set of 54 features that will serve as input to the binary classification model. In~\Cref{alg:full}, the features are computed in \cref{step:ml:get-features}.

\begin{table*}[htbp!]
    \footnotesize
	\centering
	\caption{Features for the $i^\text{th}$ constraint $A_{i\cdot} x + C_{i\cdot} v = b_i$ and a fractional solution $(x^\star,v^\star)$ of a MIP.}
			\begin{tabularx}{\textwidth}{lXr}
				\toprule
				Feature         &  Description                                              &  Count \\
				\midrule
				Right-hand side (RHS)      &  The raw value of $b_i$ and a categorical feature if it is nonzero &  2     \\
                Slack                &  The raw value of the slack variable and a categorical feature if it is nonzero &  2     \\
                Dual                 &  The value of the dual variable in the optimal LP solution &  1 \\
                Degree               &  The number of variables with nonzero coefficients in the constraint: considering all the variables, only the variables that are nonzero at $(x^\star,v^\star)$, only the variables that are zero at $(x^\star,v^\star)$, and only the variables that are at their upper bound at $(x^\star,v^\star)$ & 3 \\
                Sense                &  One-hot encoding of the sense of the constraint &  2 \\
                Stats of the coefficients &  Mean, standard deviation, minimum, and maximum of the coefficients: considering all the variables, only the variables that are nonzero at $(x^\star,v^\star)$, only the variables that are zero at $(x^\star,v^\star)$, and only the variables that are at their upper bound at $(x^\star,v^\star)$ & 16 \\
                Stats of the ratios  &  Mean, standard deviation, minimum, and maximum of the ratios between the coefficients and the RHS: considering all the variables, only the variables that are nonzero at $(x^\star,v^\star)$, only the variables that are zero at $(x^\star,v^\star)$, and only the variables that are at their upper bound at $(x^\star,v^\star)$ & 16 \\
                Euclidean distance to  $(x^\star,v^\star)$ &  From~\cite{Wess}  & 1 \\
                Relative violation &  From~\cite{Wess}  & 1 \\
                Adjusted distance to  $(x^\star,v^\star)$ &  From~\cite{Wess}  & 1 \\
                Objective function parallelism &  From~\cite{Wess} &  1 \\
                Stats of the cost vector &  Mean, standard deviation, minimum, and maximum of the cost coefficients of the variables with nonzero coefficients in the constraint &  4 \\
                                    &  We normalize the absolute value of the cost vector, and compute the number of variables on the top 1, 5, 10, and 20\% of the costs that appear with non zero coefficient in the constraint &  4 \\
                \midrule
				Total count                                &                                                                 &  54    \\
				\bottomrule
			\end{tabularx}
	\label{tab:features}
\end{table*}

As for the labels, the observation corresponding to constraint $i$ on any round is assigned a positive label if $\lvert \lambda_i \rvert >\epsilon$ for any cut in the cut pool of that round of separation, and a negative label otherwise. This gives us a dataset that can be used for the traditional binary classification task in ML.

\subsection{Solving a Reduced Separator}

We can use the output of the ML model to predict which constraints will be useful to generate MIR cuts. Then, at each iteration of the cutting loop in~\Cref{alg:full}, we compute the features for the current fractional solution on \cref{step:ml:get-features}, and use the ML model to classify each constraint on \cref{step:ml:get-constraints}. We use the predicted class for each constraint to update Appx-MIR-Sep and fix $\lambda_i=0$ for those constraints that the ML model classified in the negative class. 
We call this the \textbf{``reduced'' separator}.

\section{Computational Experiments}
\label{sec:computational-experiments}

In our experiments, we compare the ``reduced'' and ``full'' separators to evaluate how removing constraints that are predicted to be unimportant affects Appx-MIR-Sep.

\subsection{Evaluation and Computational Setup}
To evaluate the impact of a set of cuts, we use the percent integrality gap closed as defined in \Cref{sec:pool}. We implement the pipeline described above and test it on instance families derived from three base instances in the benchmark set of MIPLIB 2017~\cite{miplib}: \instance{binkar10\_1}, \instance{gen-ip054}, and \instance{neos5}. 

These three base instances are selected from 37 MIPLIB 2017 ``Benchmark'' instances for which two conditions are satisfied:
    (1) random perturbations of the objective function coefficients lead to 1,000 distinct LP relaxation optima;
    and (2) ~\Cref{alg:full} closes at least 5\% gap on average across the random variations.
We limit our experiments to only three of these 37 potential instances due to the high computational costs of running our experiments and limitations on our shared computational resources:
for each of the 1,000 instances in an instance family, 
the cut generation loop of~\Cref{alg:full} runs for up to three hours.
\Cref{tab:split} provides characteristics on the three chosen instances, which sample some of the diversity of the broader set:
\instance{binkar10\_1} and \instance{neos5} are mixed integer, whereas \instance{gen-ip054} is pure integer;
\instance{binkar10\_1} has more than 2,000 variables and 1,000 constraints, whereas the other two are smaller.

We run \Cref{alg:full} on each of the 1,000 variations from each instance family. We set the time limit to 10 minutes for each individual run of Appx-MIR-Sep, and to three hours for the entire cutting loop per instance.
The instances are then split into three disjoint parts: ``\textbf{training}'', ``\textbf{test}'', and ``\textbf{unseen}'' sets.
The training and test sets are an 80/20 split among instances in which MIR cuts perform well: we discard perturbed instances from each family that have less than 5\% gap closed at the end of the cutting loop.
\Cref{tab:split} shows the number of instances in each of the sets.

\begin{table}[t]
	\centering
	\caption{Descriptive statistics for the three instance families.}
	\label{tab:split}
	\begin{tabular}{ l@{\hskip .1in}
        *{2}{S[table-format=4.0]}
        *{4}{S[table-format=3.0]}
    }
		\toprule
		Family      & {Constraints} & {Continuous vars.} &  {Integer vars.}  &  {Training} & {Test} & {Unseen}  \\
		\midrule
		\instance{binkar10\_1} & 1026        & 2128             & 170             & 478        & 120 & 402 \\
		\instance{gen-ip054} & 27          & 0                & 30              & 416        & 104 & 480 \\
		\instance{neos5}       & 63          & 10               & 53              & 611        & 153 & 236 \\
		\bottomrule
	\end{tabular}
\end{table}

\mbox{}\\[-1ex]\noindent\textbf{Training set:}
For each instance family, we use the corresponding training data to fit a gradient-boosted tree ensemble using \texttt{sklearn}'s \cite{scikit-learn} implementation with default hyperparameters except for the maximum tree depth, which we switch from 3 to 5. \Cref{tab:acc} shows that the trained models perform fairly well as measured by their accuracy, precision (fraction of constraints predicted to be useful that are labeled as useful), and recall (fraction of constraints labeled as useful that are predicted to be useful). %
For example, for instance \texttt{binark10\_1}, 80\% of the constraints that the classifier selects to keep have a label of ``useful'', and 93\% of the constraints with a label of ``useful'' are correctly classified.

\mbox{}\\[-1ex]\noindent\textbf{Test set:}
We next use the trained ML models to reduce the Appx-MIR-Sep problem in each round of separation. Since we know \latin{a priori} that the reduced separator closes more than 5\% gap in these variations, this comparison favors the full separator. 

\begin{table}[htbp!]
	\centering
	\caption{Trained classifier performance on the training and test sets.}
	\label{tab:acc}
	\begin{tabular}{ l 
        @{\hskip .1in}c @{\hskip .1in}c @{\hskip .1in}c
        @{\hskip .1in}c @{\hskip .1in}c @{\hskip .1in}c 
    }
		\toprule
        & \multicolumn{3}{c}{Training} 
        & \multicolumn{3}{c}{Test} \\
            \cmidrule(r){2-4} \cmidrule{5-7}
		Family      & Accuracy & Precision & Recall 
            & Accuracy & Precision & Recall 
        \\
		\midrule
		\instance{binkar10\_1} & 0.819    & 0.786     & 0.928  
            & 0.823    & 0.794     & 0.926  
        \\
		\instance{gen-ip054} & 0.878    & 0.839     & 0.778  
            & 0.820    & 0.674     & 0.846  
        \\
		\instance{neos5}       & 0.940    & 0.776     & 0.820  
            & 0.936    & 0.755     & 0.833  
        \\
		\bottomrule
	\end{tabular}
\end{table}

\mbox{}\\[-1ex]\noindent\textbf{Unseen set:}
The last set is of instances from each family not considered in the training and test sets, i.e., objective functions that lead to an initial fractional optimal solution from which the eventual gap closed is less than $5\%$ at the end of the cutting loop.
Effectively, we investigate if this negative performance is inherent to these instances, i.e., we cannot find good MIR cuts from these starting points, or if the learned ``useful'' constraint set produces better MIR cuts from Appx-MIR-Sep.

\subsection{Computational Results}
\label{sec:computational-results}

\Cref{fig:binkar,fig:gen,fig:neos} show the effect of the reduced separator for each instance family.
The three panels correspond to the training set on the left, test set in the middle, and unseen set on the right.
The horizontal axis is the round of separation.
The bottom subplots show the number of instances for which the cutting loop does not terminate by that round.
In the top subplots, the vertical axis is percent integrality gap closed.
The solid lines show the average percent gap closed calculated over only the instances that require at least that many rounds of the cutting plane loop.
This set of surviving instances may differ across the methods, and the behavior is nonmonotone since an instance with a high gap closed no longer contributes to the average after its last round of cuts.
The error bars show the standard deviation of the percent gap closed over the surviving instances.
The dashed lines show the average percent gap closed across all instances in the family,
where the gap closed for an instance that terminates at an earlier round remains constant for all subsequent rounds.

We will also refer to
\Cref{tab:data} for averaged statistics for each instance family, dataset, and method for
    (1) number of cutting plane rounds, 
    (2) number of generated cuts,
    and
    (3) time for cut generation.
Generally, the reduced separator uses fewer rounds, generates a comparable number of cuts, and requires more time.
We do not focus on cut generation time, as it can be easily skewed: an instance that does no rounds of cuts at all would have a lowest-possible value of zero on that metric.

{
\sisetup{
    table-alignment-mode = format,
    table-number-alignment = center,
    table-format = 2.2,
    table-auto-round,
}
\begin{table}
\centering
\caption{Average of number of rounds, number of generated cuts, and seconds to solve Appx-MIR-Sep for each dataset.}
\label{tab:data}

\begin{tabular}{
  l
  l
  S[table-format=2.2]
  S[table-format=2.2]
  S[table-format=2.2]
  S[table-format=2.2]
  S[table-format=3.2]
  S[table-format=3.2]
}
\toprule
\multicolumn{2}{c}{} & \multicolumn{2}{c}{Rounds (\#)} & \multicolumn{2}{c}{Cuts (\#)} & \multicolumn{2}{c}{Time per round (s)} \\
\cmidrule(lr){3-4} \cmidrule(lr){5-6} \cmidrule(lr){7-8}
 {Family} & {Set}  & {Full} & {Reduced} & {Full} & {Reduced} & {Full} & {Reduced} \\
\midrule
\instance{binkar10\_1}
    & Training & 24.008368 & 10.506276 & 16.024922 & 16.148945 & 307.242168 & 496.794729 \\
    & Test & 24.475000 & 11.133333 & 16.403473 & 16.308383 & 296.655228 & 501.845475 \\
    & Unseen & 11.873134 & 10.646766 & 15.421328 & 16.428037 & 231.241322 & 492.733497 \\
\addlinespace
\instance{gen-ip054}
    & Training & 4.920673 & 4.466346 & 16.138740 & 15.832078 & 279.060663 & 270.323810 \\
    & Test & 4.980769 & 4.298077 & 15.789575 & 15.548098 & 279.619658 & 289.077633 \\
    & Unseen & 2.962500 & 3.475000 & 14.397328 & 14.661271 & 172.848718 & 194.789595 \\
\addlinespace
\instance{neos5}
    & Training & 2.787234 & 1.931260 & 5.266001 & 4.610169 & 107.153939 & 226.923759 \\
    & Test & 2.588235 & 1.888889 & 5.250000 & 4.899654 & 103.110461 & 255.446518 \\
    & Unseen & 1.072340 & 1.421277 & 4.932540 & 4.613772 & 0.018496 & 235.401920 \\
\bottomrule
\end{tabular}
\end{table}
} %

\subsubsection{Instance Family \instance{binkar10\_1}}
\label{sec:binkar10_1}
\begin{figure}
    \centering
	\includegraphics[width=\figurewidth]{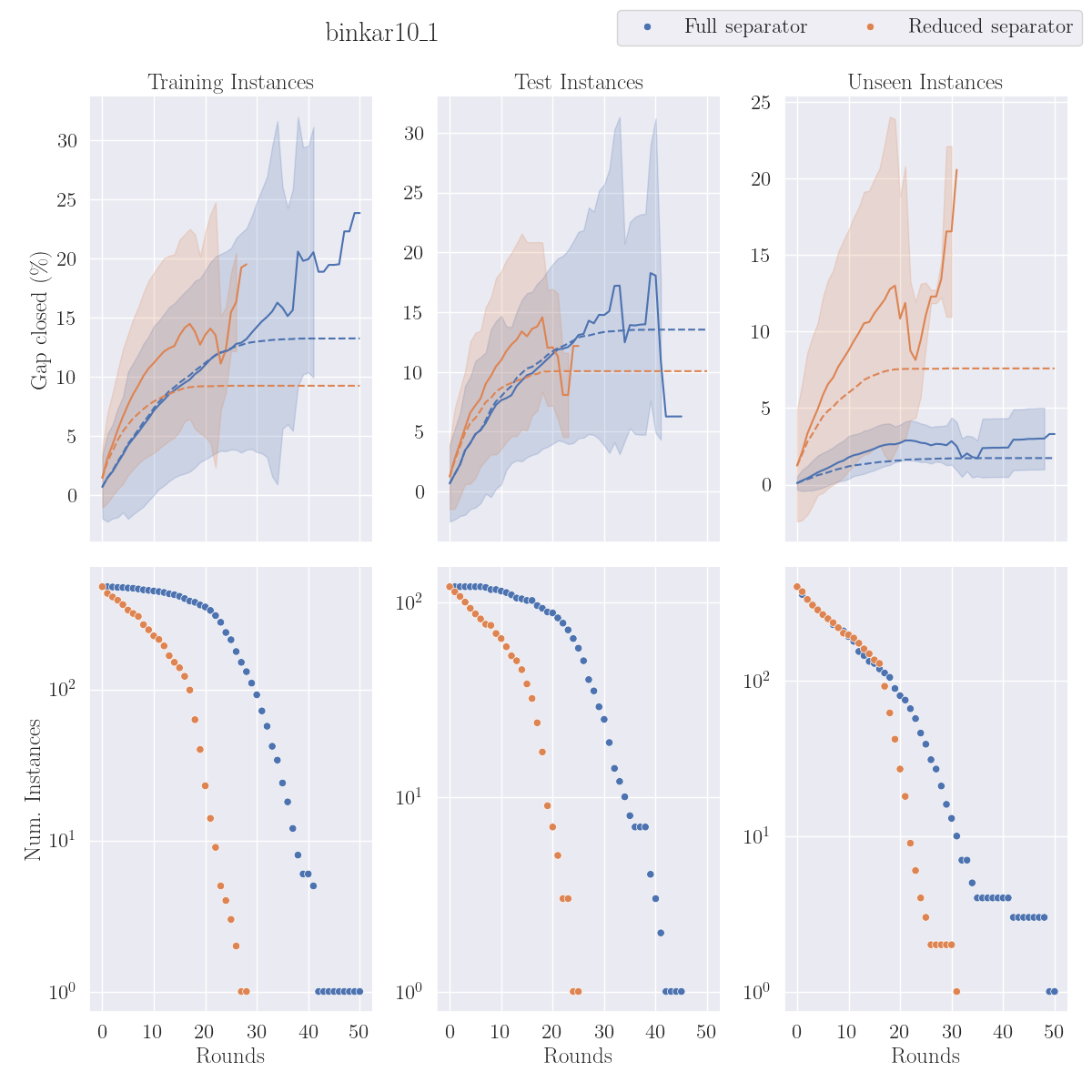}
	\caption{%
        Effect of reduced separator on instance family \instance{binkar10\_1}. %
    } \label{fig:binkar}
\end{figure}

\Cref{fig:binkar} shows the results for instance family \instance{binkar10\_1}.
We first discuss the training and test sets.
The reduced separator tends to close a larger percent of the integrality gap
for the instances that survive to a given round 
(these may not be the same subset across methods),
However, as seen in the bottom subplots and in \Cref{tab:data}, the cutting loop with the reduced separator tends to terminate earlier.
This means that the reduced separator terminates for many instances at an early round with a small gap closed.
The reduced separator may close more gap by that round for those instances,
but the full generator continues for more rounds and eventually tends to provide better bounds.
This is seen with the dashed lines in the top subplots.

The reduced separator is about 30\% smaller than the full separator for this family, for both training and test sets.
However, from \Cref{tab:data}, we see that cut generation time is larger on average.
Hence, the smaller Appx-MIR-Sep requires more time to solve to optimality.
The reduced separator generates around the same number of cuts in total as the full separator,
but it obtains these over fewer rounds.
Thus, the reduced separator finds more cuts in each round, on average, which may explain the observed larger gap closed (for the surviving instances per round).

\Cref{tab:data} does not state the cumulative time spent on cut generation per method.
With respect to that metric, for the training and test sets, since the reduced generator requires significantly fewer rounds, it also takes 25--30\% less time on average, even though each round is more expensive.

For the unseen set, the reduced separator, using 71\% of the constraints, closes more gap than the full separator for all variations, with (slightly) fewer rounds of the cutting loop,
but on average Appx-MIR-Sep running time (both per round and, for this set, cumulative) is higher.
This additional cost is related to (and justified by) the higher ``success'' (in terms of cut quality) solving Appx-MIR-Sep for these instances.

\subsubsection{Instance Family \instance{gen-ip054}}
\label{sec:gen-ip054}

\begin{figure}
    \centering
	\includegraphics[width=\figurewidth]{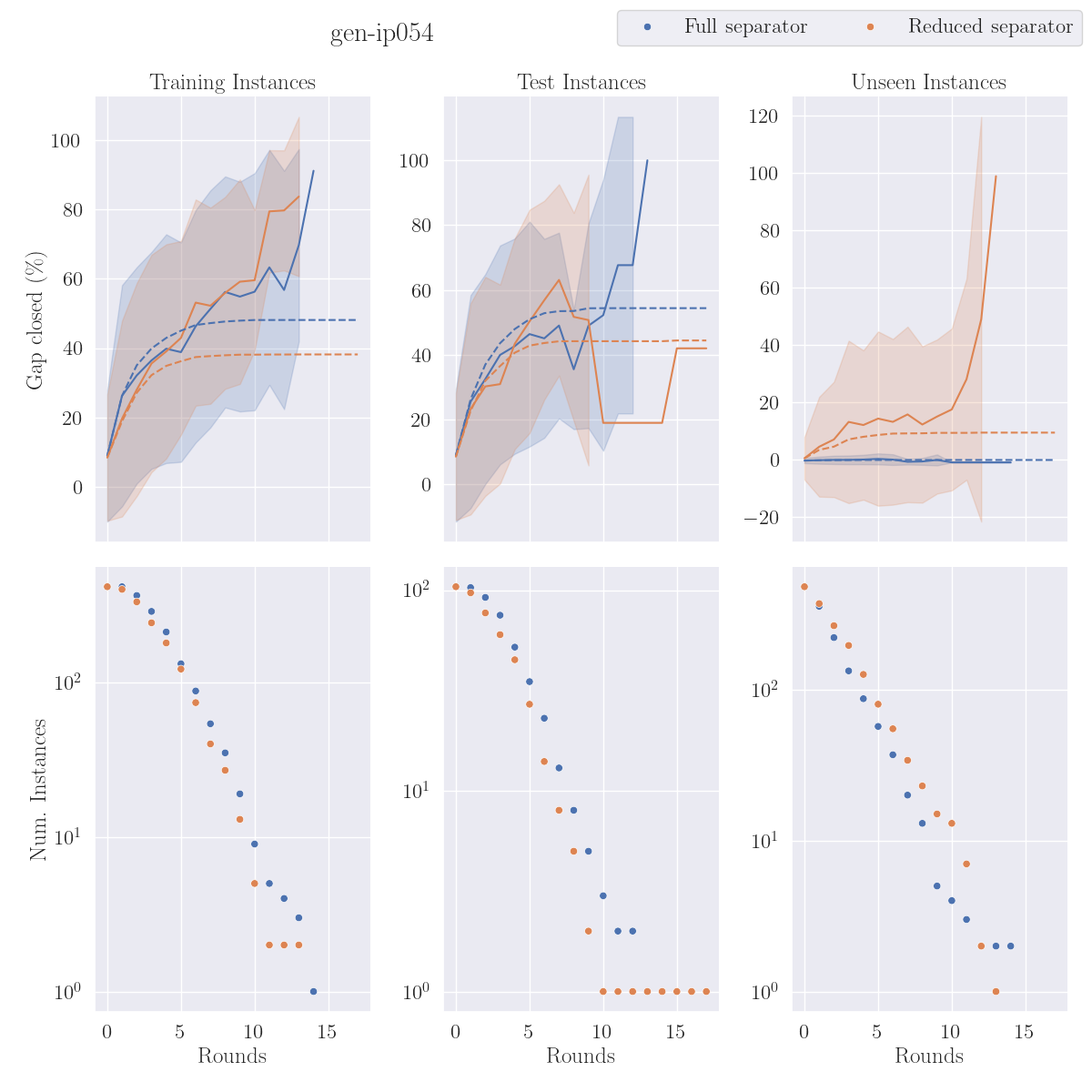}
	\caption{%
        Effect of reduced separator on instance family \instance{gen-ip054}. %
    } \label{fig:gen}
\end{figure}

\Cref{fig:gen} shows the results for instance family \instance{gen-ip054}.
We again first analyze the training and test sets.
Both generators have a similar initial progression of average percent gap closed for surviving instances in these sets. 
The cutting loop for the reduced separator terminates after slightly fewer iterations in most cases.
However, on the instances that do not terminate early, the full separator closes more gap. 
Around 30\% of the constraints are labeled ``unimportant'' by the classifier for both training and test sets.
In contrast to the situation with \instance{binkar10\_1}, for this family we observe in \Cref{tab:data} that there is a slight decrease in time to solve Appx-MIR-Sep for the training set,
and a slight increase for the test set,
but the values are comparable.

For the unseen set, the reduced separator, using 75\% of the constraints, again outperforms the full separator in terms of gap closed (both on average and per round) and notably closes 100\% of the gap for one instance.
The reduced separator also performs more cutting plane rounds on average than the full generator,
and each round requires more time,
leading to a 30\% increase in cumulative time on cut generation.

\subsubsection{Instance Family \instance{neos5}}
\label{sec:neos5}

\begin{figure}
    \centering
	\includegraphics[width=\figurewidth]{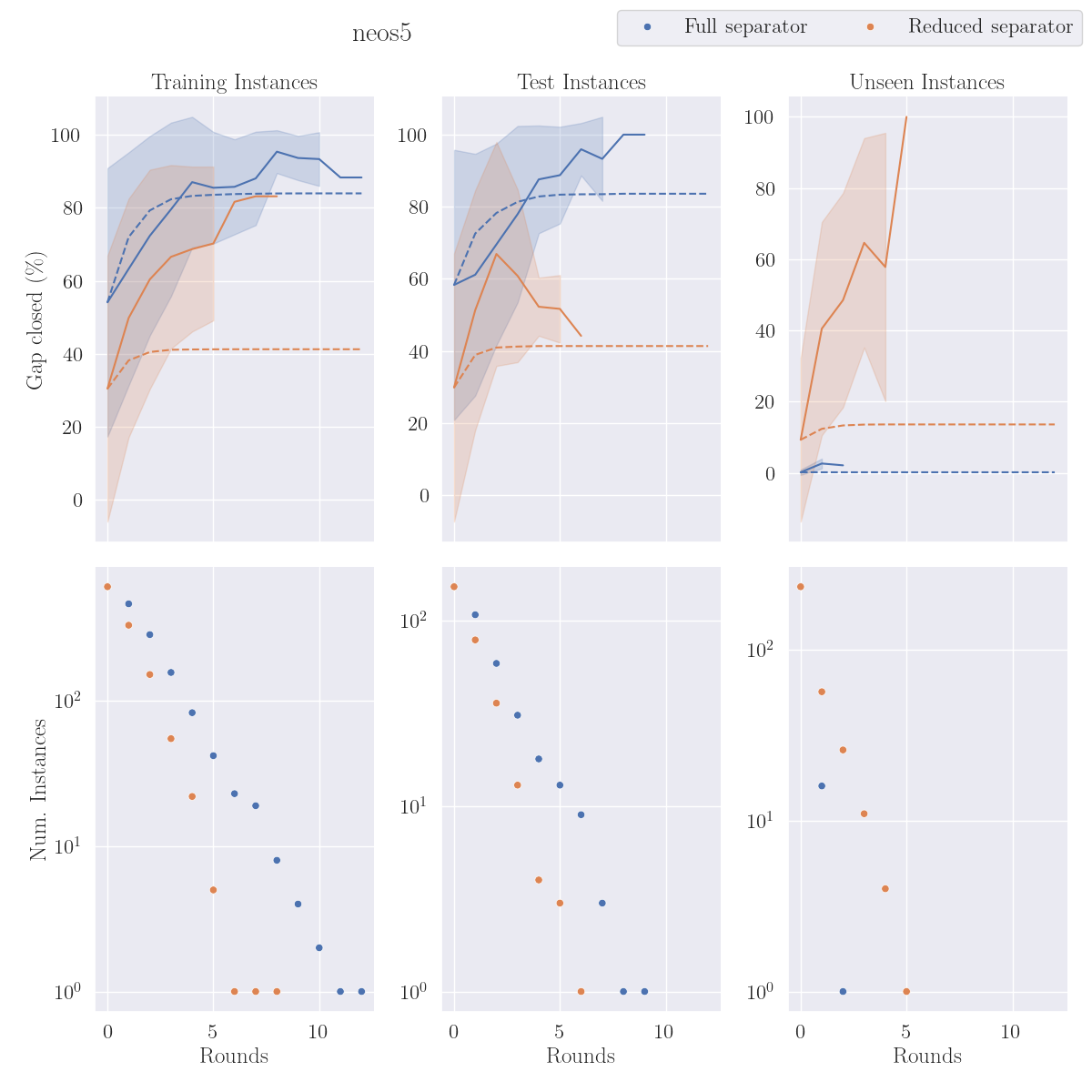}
	\caption{%
        Effect of reduced separator on instance family \instance{neos5}. %
    } \label{fig:neos}
\end{figure}

Finally, \Cref{fig:neos} shows our results for instance family \instance{neos5}.
We observe that the full separator closes more gap than the reduced separator.
On average, for the training and test sets, only 15\% of the constraints are labeled ``useful'',
but, from \Cref{tab:data}, the average cut generation time (per round), i.e., the average time to prove optimality of Appx-MIR-Sep, is still longer for the reduced separator.
For the unseen set,  the reduced separator uses 16\% of the constraints and consistently outperforms the full separator, closing 100\% of the gap for one instance.

\section{Conclusion}

\begin{itemize}[leftmargin=*]
    \item 
        [\textbf{Takeaway 1:}]
    
    \Cref{fig:binkar,fig:gen,fig:neos} show that the reduced separator tracks the full separator reasonably well on both training and test instances. Importantly, these instances are ones for which we know that the full separator closes at least 5\% of the gap. Since the ML model is trained to imitate the cuts produced by the full separator on the training instances, it is not surprising that the performance of the reduced separator is typically upper bounded by that of the full separator.\\

    \item
        [\textbf{Takeaway 2:}]
    The results shown for the unseen set in \Cref{fig:binkar,fig:gen,fig:neos} are especially promising, as they show the potential for incorporating ML into optimization-based separators. When optimizing purely for violation, the total gap closed by the full separator is negligible, but when using the information learned by the ML models, the (reduced) separator generates cuts that close much more gap, the true goal in cut generation.
\end{itemize}

We have explored the potential of a rather simple binary classification formulation for the problem of reducing the set of constraints that are considered by an optimization-based cut separator. While more sophisticated ML models such as graph neural networks or transformers could be used instead of a gradient-boosted tree ensemble, we have chosen the latter for its simplicity in terms of fitting and amount of training data it requires. We are yet to achieve a significant improvement in the running time of the reduced MIR separator as compared to the full separator. We believe this to be an important step towards operationalizing this hybrid approach. Although we have considered only three MIPLIB 2017-derived instance datasets, we note that generating the training data and evaluating both the full and reduced separators on the training and test instances required hundreds of CPU days in computation. Expanding the experiments to more instance families both from MIPLIB 2017 and other more structured problem classes is of immediate interest.

{
\small
\mbox{}\\\noindent\textbf{Acknowledgments}
Aleksandr M. Kazachkov acknowledges the following funding support: This material is based upon work supported by the Air Force Office of Scientific Research under award number FA9550-23-1-0340.
Elias Khalil acknowledges funding from the Natural Sciences and Engineering Research Council of Canada - Discover Grant Program and from the SCALE AI Research Chair Program.
This research was enabled in part by support provided by Compute Ontario (\url{https://www.computeontario.ca}) and the Digital Research Alliance of Canada (\url{https://alliancecan.ca}).

\mbox{}\\\noindent\textbf{Competing Interests}
The authors have no relevant financial or non-financial interests to disclose.

\mbox{}\\\noindent\textbf{Data Availability Statement}
All data supporting the findings of this study are available via the University of Toronto Dataverse at \url{https://doi.org/10.5683/SP3/KBVV6C}. Code for this paper is at \url{https://github.com/khalil-research/LearnMIR/}.
}

\bibliographystyle{plainnat}
\bibliography{mybibliography}%

\end{document}